\documentclass[a4paper,12pt]{article}
\usepackage{amssymb,amsmath,amsthm,times}
\newtheorem{thm}{Theorem}[section]

\newtheorem{lem}[thm]{Lemma}

\def\ac#1#2{\binom{#1 #2\theta -1}{ #1}}
\def\r1{|1-z|}

 \title{   Ces\`{a}ro summation and multiplicative functions on a symmetric group}

 \author {Vytas Zacharovas}

\begin{document}

\maketitle

\begin{abstract}

We investigate the summability in sense of Ces\`{a}ro and its
applications to investigation of the mean values of multiplicative
functions on permutations.
\end{abstract}

\section{Results}
Let $S_n$ be the symmetric group. Each element
 $\sigma \in S_n$ can be
 decomposed into a product of independent
 cycles.
 \[
 \sigma=\kappa_1\kappa_2...\kappa_{\omega}
 \]
this decomposition is unique up to the order of the multiplicands.
We will call a function $f: S_n\to C$ multiplicative if
$f(\sigma)=f(\kappa_1)f(\kappa_2)...f(\kappa_n)$. In what follows we
will assume that the value of $f$ on cycles depends only on the
length of cycle, that is $f(\kappa)=\hat f(|\kappa |)$, where
$|\kappa|$ - the order of cycle $\kappa$.
 Let
 $m_k(\sigma)$ be equal to the number of cycles in the decomposition  $\sigma$
 whose order is equal to $k$. Then obviously $m_1(\sigma)+
2m_2(\sigma)+...+nm_n(\sigma)=n$. Thus
  $n$ complex number $\hat f(1),\hat f(2),...,\hat f(n)$ completely
  determine the value of function $f$ on any permutation $\sigma \in S_n$
$$
f(\sigma)=\hat f(1)^{m_1(\sigma)}\hat f(2)^{m_2(\sigma)}...
\hat f(n)^{m_n(\sigma)}.
$$
On the group $S_n$ we will define the so called Ewens's measure
$\nu_{n,\theta}$ by means of formula
\[
\nu_{n,\theta}(\sigma)={{\theta^{k(\sigma)}}\over {\theta_{(n)}}},
\]
where $k(\sigma)=m_1(\sigma)+m_2(\sigma)+...+m_n(\sigma)$, and
$\theta_{(n)}=\theta(\theta+1)...(\theta+n-1)$.

We will investigate the mean values of multiplicative functions with
respect to Ewens measure
\[
M_n(f)=\sum_{\sigma \in S_n}f(\sigma)\nu_{n,\theta}(\sigma).
\]
Since the number of $\sigma$ such that  $m_k(\sigma)=s_k$  is equal
to $ n!\prod_{j=1}^{n} {1 \over {s_j!j^{s_j}}}, $ therefore
\[
\nu_{n,\theta}(m_1(\sigma)=s_1,...,m_n(\sigma)=s_n)={{n!}\over{\theta_{(n)}}}
\prod_{j=1}^{n} {\left({\theta \over j}\right)}^{s_j}{1\over{s_j!}}.
\]
Hence
\[
M_n(f)={\ac{n}{+}}^{-1} \sum_{k_1+2k_2+...+nk_n=n} {\prod_{j=1}^{n}}
{     {     \left(     {{\theta \hat f(j)}\over j}     \right)
}^{k_j}    } {1 \over {k_j!}}.
\]
It is easy to see that $M_n(f)$ is equal to ${\ac{n}{+}}^{-1}N_n$,
where $N_n$ is defined by means of relation
\[
F(z)=\exp \left\{\theta \sum_{j=1}^{\infty} {{\hat f(j) \over j}z^j}
\right\}=\sum_{m=0}^{\infty}N_mz^m.
\]
Since the numbers  $\hat f(j)$ with $j>n$  do not influence the
value of the coefficient of $z^n$ therefore we will assume that
$\hat f(j)=1$ for $j>n$. Therefore
\[
F(z)=\exp \left\{\theta \sum_{j=1}^{\infty} {{\hat f(j) \over j}z^j}
\right\}=\sum_{j=0}^{\infty}N_jz^j= {{\exp\{ \theta L_n(z) \}}\over
{(1-z)^{\theta}}} ,
\]
here and in what follows $L_n(z)=\sum_{j=1}^{n}{{\hat f(j)-1}\over
j}z^j$, and $L_0(z)=0$.

The function $F(z)$ is the product of two functions $\exp \{ \theta
L_n(z)\}=\sum_{k=0}^{\infty}m_kz^k$ and ${1\over
{(1-z)^{\theta}}}=\sum_{k=0}^{\infty}\ac{n}{+}z^n,$ therefore
\begin{equation}
\label{M_n} M_n(f)={\ac{n}{+}}^{-1}\sum_{j=0}^nm_j\ac{n-j}{+}.
\end{equation}

We will estimate the sum on the right hand side of the equation
(\ref{M_n}) by means of the following theorem.

\begin{thm}
\label{T_main_thm} Let  $f(z)=\sum_{m=0}^{\infty}a_mz^m$ be analytic
for  $|z|<1$. Let us denote
$$
S_{\theta}(f;n)=\sum_{k=1}^nka_k\ac{n-k}{+},
$$
then for fixed  $\theta >0$ we have
\begin{multline*}
{1\over {\ac{n}{+}}}\sum_{k=0}^na_k\ac{n-k}{+}-f(e^{-1/n})-
{{S_{\theta}(f;n)}\over {n\ac{n}{+}}}
\\
\ll{1\over {n}}\sum_{j=1}^{\infty}
{{|S_{\theta}(f;j)|}\over {j^{\theta}}}e^{-j/n}
+{1\over {n^\theta}}
\sum_{j=n}^{\infty}{{|S_{\theta}(f;j)|}\over {j}}e^{-j/n}.
\end{multline*}
The constant in the symbol $\ll$ depends only on  $\theta$.
\end{thm}

The sum on the righthand side of   (\ref{M_n}), is called Ces\`{a}ro
mean with parameter $p=\theta -1$. If for a given formal series
$\sum_{k=0}^{\infty}a_k$ the Cesaro means with parameter $p$
converge to some number $A$, then we say that
$\sum_{j=0}^{\infty}a_k$ is  $(C,p)$ summable and its Ces\`{a}ro sum
is $A$ and write $(C,p)\sum_{j=0}^{\infty}a_j=A$.

From Theorem \ref{T_main_thm} we can deduce the following result,
which is probably already known.

\begin{thm}
\label{T_tauberian} Suppose $p>-1$. A series
$\sum_{k=0}^{\infty}a_k$ is $(C,p)$ with  summable and it's $(C,p)$
sum is equal to $A$ if and only if
\begin{equation}
\label{first} \lim_{x \to 1-0}\sum_{k=0}^{\infty}a_kx^k=A,
\end{equation}

\begin{equation}
\label{second} \lim_{n\to \infty }{{S_{p+1}(f;n)}\over {n^{p+1}}}=0,
\end{equation}
where $f(x)=\sum_{j=0}^{\infty}a_jx^j$.
\end{thm}

In the case when $\theta =1$ Theorem \ref{T_tauberian} becomes the
classical theorem of Tauber (see. e.g.
\cite{tenenbaum},\cite{hardy}). For this special case the proof of
Theorem \ref{T_main_thm} can be obtained by modifying the proof of
Tauber's theorem. Let us define
$$
\mu_n(p)={\left( {1\over n}\sum_{k=1}^n|\hat f(k)-1|^p
\right)}^{1/p}.
$$
Applying Theorem \ref{T_main_thm} we can easily prove the following
result
\begin{thm}
\label{T_mean_value}
 Пусть $p>\max\left\{ 1,{1\over {\theta}} \right\}$ и $|\hat
f(j)| \leqslant 1$, тогда
\[
M_n(f)=\exp \left\{ \theta \sum_{k=1}^{n}{{\hat f(k)-1}\over k}
\right\} +O\bigl( \mu_n(p) \bigr),
\]
here the constant in symbol $O(..)$ depends only on $\theta$ and
$p$.
\end{thm}
The variants of  Theorem \ref{T_mean_value} with less precise
estimate of the remainder term have been proved in
\cite{manst_berry-esseen},\cite{manst_tauber},\cite{manst_decomposable}
and \cite{manst_additive}.

\section{ Proofs}

\begin{lem}
\label{L_flaj_odl} Let $f(z)=\sum_{m=0}^{\infty}f_nz^n$ be analytic
function in the region $\Delta (\phi,\eta )=\{ z |\quad |z|<1+\eta,
\quad |\arg(z-1)|>\phi \}$, where
 $\eta>0$ and $0<\phi <\pi /2$. If
$$
|f(z)|\leq K_1 |1-z|^{\alpha_1}+K_2 |1-z|^{\alpha_2},
$$
for $z\in \Delta (\phi,\eta )$ then there exists such a constant $c=
c(\alpha_1,\alpha_2,\eta,\phi)$ which is independent of $K_1, K_2$
and such that
$$
|f_n|\leq c(K_1n^{-\alpha_1-1}+K_2n^{-\alpha_2-1}).
$$
\end{lem}
\begin{proof} The same as of Theorem 1 of \cite{flajolet_odlyzko}.

\end{proof}

Let us denote
$$
c_{m,j}=\sum_{s=0}^m{{\ac{m-s}{+}\ac{s}{-}}\over{s+j}},
$$
for $j\geq 1$. Then the generating function of  $c_{m,j}$ will have
the form
$$
F_j(z)=\sum_{m=0}^{\infty}c_{m,j}z^m={1\over{(1-z)^{\theta}}}\int_{0}^1
(1-xz)^{\theta}x^{j-1}dx.
$$

\begin{lem}
\label{L_estimate_c_m_j} We have the following estimates for
$c_{m,j}$ :
$$
\leqno (i)\quad 0\leq c_{m,j} \leq {\theta \over{j^2}} e^{\theta
m/j},\quad m\geq 1, \quad c_{0,j}={1\over j} ;
$$
$$
\leqno (ii)\quad
c_{m,j}=\ac{m}{+}\int_0^1(1-y)^{\theta}y^{j-1}dy+O\left(
{m^{\theta-2}\over {j^\theta}}+ {1\over {m^2}}\right).
$$
\end{lem}
\begin{proof} Differentiating $F_j(z)$ we obtain
$$
zF_j'(z)={{\theta z F_j(z)}\over {1-z}}+1-j F_j(z).
$$
Expanding both sides of the above equation into Taylor series and
equating the coefficients of the same powers $z^m$ we obtain
$$
c_{m,j}={{\theta}\over{m+j}}\sum_{s=0}^{m-1}c_{s,j},\quad m\geq 1
$$
and $c_{0,j}={1\over j}$. This recurrent relation implies that
$$
0<c_{m,j}\leq {{\theta}\over{j}}\sum_{s=0}^{m-1}c_{s,j},\quad m\geq 1.
$$
Then
$$
0\leq c_{m,j}\leq b_{m,j},
$$
where $b_{m,j}$ is solution of the recurrent equation
$$
b_{m,j}= {{\theta}\over{j}}\sum_{s=0}^{m-1}b_{s,j},\quad m\geq 1 ,
$$
with initial condition $b_{0,j}={1\over j}$. It is easy to check
that
$$
b_{m,j}={\theta \over {j^2}}{\left( 1+{\theta \over j} \right)}^{m-1},\quad
m\geq 1.
$$
Therefore applying inequality $1+x\leq e^x$ we obtain the estimate
$(i)$
$$
c_{m,j}\leq  {\theta \over {j^2}}{\left( 1+{\theta \over j}
\right)}^{m-1} \leq  {\theta \over {j^2}}e^{\theta m/j} ,\quad m\geq
1.
$$

In order to prove estimate $(ii)$ we will use Lemma \ref{L_flaj_odl}
with $\eta =1/2$ and
 $\phi=\pi /4$.
We can represent $F_j(z)$  as a sum of two functions
$$
F_j(z)={1\over {(1-z)^{\theta}}}\int_{0}^1(1-x)^{\theta}x^{j-1}dx+G_j(z).
$$

Let $z\in \Delta (1/2,\pi /4)$, $|z-1|<1/2$. Then
\[
\begin{split}
\int_{0}^1&(1-zy)^{\theta}y^{j-1}dy-\int_{0}^1(1-y)^{\theta}y^{j-1}dy
\\
&= \int_0^{1-\r1}(1-y)^{\theta}\left( {\left( 1-y{{z-1}\over {1-y}}
\right)}^{\theta}-1 \right) y^{j-1}dy+
\\
&\quad+\int_{1-\r1}^1 \left( (1-y-y(z-1))^{\theta}-(1-y)^{\theta}
\right) y^{j-1} dy
\\
&\ll \int_0^{1-\r1 }(1-y)^{\theta}{{y^j\r1}\over {1-y}}dy
+\int_{1-\r1}^1y^{j-1}{\r1}^{\theta}dy
\\
&\ll\r1 \int_0^1(1-y)^{\theta-1}y^{j-1}dy +
{\r1}^{\theta+1}
\\
&\ll{{|1-z|}\over {j^{\theta}}}+|1-z|^{\theta+1}.
\end{split}
\]

It is easy to see that the obtained estimate holds in the whole
region $\Delta(\eta,\phi)$. Therefor for $z\in \Delta (\eta,\psi)$
\begin{equation}
\label{G_j}\begin{split} G_j(z)&={1\over
{(1-z)^{\theta}}}\int_0^1\left( (1-yz)^{\theta}-(1-y)^{\theta}
\right) y^{j-1}dy
\\
&\ll|1-z|+{{|1-z|^{1-\theta}}\over {j^{\theta}}}.
\end{split}
\end{equation}
Applying Lemma \ref{L_flaj_odl} with $f(z)=G_j(z)$ and taking into
account (\ref{G_j}) we obtain estimate $(ii)$.

The lemma is proved.
 \end{proof}

\begin{proof}[ Proof of Theorem \ref{T_main_thm}] Since
$$
\sum_{k=1}^{\infty}S_{\theta}(f;k)z^k={{zf'(z)}\over
{(1-z)^{\theta}}},
$$
then
$$
na_n=\sum_{k=1}^nS_{\theta}(f;k)\ac{n-k}{-},\quad n\geq 1.
$$
Therefore
\[
\begin{split}
R_n&:=\sum_{k=0}^na_k\ac{n-k}{+}-f(e^{-1/n})\ac{n}{+}
\\
&=\sum_{k=1}^n\ac{n-k}{+}{1\over
k}\sum_{j=1}^kS_{\theta}(f;j)\ac{k-j}{-}
\\
&\quad-\ac{n}{+}\sum_{k=1}^{\infty}e^{-k/n}{1\over
k}\sum_{j=1}^kS_{\theta}(f;j) \ac{k-j}{-}
\\
&=\sum_{j=1}^{n}S_{\theta}(f;j)c_{n-j,j}
-\ac{n}{+}\sum_{j=1}^{\infty}S_{\theta}(f;j)\sum_{k=j}^{\infty}
{{\ac{k-j}{-}e^{-k/n}}\over k}.
\end{split}
\]

Suppose $j>n/2$, then
\[
\begin{split}
\sum_{k=j}^{\infty} {{\ac{k-j}{-}e^{-k/n}}\over
k}&=\sum_{s=0}^{\infty}{{\ac{s}{-}e^{-{{j+s}\over n}}}\over
{j+s}}=\int_0^{e^{-1/n}}(1-x)^{\theta}x^{j-1}dx
\\
&=\int_{1/n}^{\infty}(1-e^{-y})^{\theta}e^{-jy}dy\leq
\int_{1/n}^{\infty}y^{\theta}e^{-jy}dy
\\
&\ll{{e^{-j/n}} \over {jn^{\theta }}}.
\end{split}
\]
Applying the obtained estimate and Lemma \ref{L_estimate_c_m_j} we
obtain
\[
\begin{split}
R_n&-{{S_{\theta}(f;n)}\over n}
\\
&\ll \sum_{j\leq
n/2}|S_{\theta}(f;j)|\left| c_{n-j,j}-\ac{n}{+}\sum_{k=j}^{\infty}
{{\ac{k-j}{-}e^{-k/n}}\over k}\right|
\\
&\quad+\sum_{n/2<j<n}|S_{\theta}(f;j)|c_{n-j,j}
+\ac{n}{+}\sum_{j>n/2}|S_{\theta}(f;j)| \left| \sum_{k=j}^{\infty}
{{\ac{k-j}{-}e^{-k/n}}\over k}\right|
\\
&\ll\ac{n}{+} \sum_{j\leq
n/2}S_{\theta}(f;j)\int_0^1(1-y)^{\theta}y^{j-1}dy \left|
{{\ac{n-j}{+}}\over {\ac{n}{+}}}-1 \right|
\\
&\quad+\ac{n}{+}\left( {1\over n}\sum_{j\leq
n/2}{{|S_{\theta}(f;j)|}\over {j^{\theta}}} +   {1\over
n}\sum_{n/2<j< n}{{|S_{\theta}(f;j)|}\over {n^\theta }}+    {1\over
{n^{\theta}}}\sum_{j>n}{{|S_{\theta}(f;j)|}\over {j}}e^{-j/n}\right)
\\
&\ll\ac{n}{+}{1\over n} \sum_{j=1}^{\infty}{{|S_{\theta}(f;j)|}\over
{j^{\theta}}}e^{-j/n}+ \ac{n}{+}{1\over {n^\theta}}
\sum_{j=n}^{\infty}{{|S_{\theta}(f;j)|}\over {j}}e^{-j/n},
\end{split}
\]
here we have used the fact that $\ac{n}{+}={{n^{\theta -1}}\over
{\Gamma (\theta)}}\left( 1+O \left( {1\over n} \right) \right)$.

The theorem is proved.
\end{proof}
\begin{proof}[ Proof of Theorem \ref{T_tauberian}] The sufficiency of conditions (\ref{first}) and (\ref{second})
follows immediately from  Theorem \ref{T_main_thm}. The fact that
Ces\`{a}ro summability implies (\ref{first}) and (\ref{second}) is
proved in \cite{hardy}.
\end{proof}
\begin{proof}[ Proof of Theorem \ref{T_mean_value}] Let us apply Theorem \ref{T_main_thm}
with $f(z)=\exp \{ L_n(z) \}=\sum_{k=0}^{\infty}m_kz^k$. Then
$$
\sum_{k=1}^{\infty}S_{\theta}(f;k)z^k={{zf'(z)}\over {(1-z)^{\theta}}}=
F(z)\theta \sum_{k=1}^n \bigl(\hat f(k)-1 \bigr) z^k,
$$
therefore
$$
S_{\theta}(f;m)=\sum_{k=1}^m\bigl( \hat f(k)-1 \bigr) N_{m-k}.
$$
Since $|N_k|\leq \ac{k}{+}$ then applying Cauchy inequality with
parameters ${1\over p}+{1\over q}=1$ we obtain
\[
\begin{split}
|S_{\theta}(f;m)| &\ll \left( \sum_{k=1}^m|\hat f(k)-1|^p
\right)^{1/p} \left( \sum_{k=1}^m k^{(\theta -1)q} \right)^{1/q}
\\
&\ll m^{\theta} \left( {n \over m}\right)^{1/p}\mu_n(p).
\end{split}
\]

Applying Theorem \ref{T_main_thm} and using the estimate $\exp\{
\theta L_n(e^{-1/n}) \}= \exp \{ \theta L_n(1) \}\bigl(
1+O(\mu_n(p))\bigr)$ we get
\[
\begin{split}
{{N_n}\over {\ac{n}{+}}}-\exp \left\{ \theta \sum_{k=1}^{n}{{\hat
f(k)-1}\over k} \right\} &\ll\mu_n(p)+ {{\mu_n(p)}\over n}
\sum_{m=1}^{\infty}
 \left( {n \over m}\right)^{1/p}e^{-m/n}
\\
&\quad+\mu_n(p){1\over {n^{\theta}}}\sum_{m=n}^{\infty}m^{\theta -1}
 \left( {n \over m}\right)^{1/p}e^{-m/n}
 \\
&\ll\mu_n(p).
\end{split}
\]

The theorem is proved.
\end{proof}

\end{document}